\documentclass[reqno]{amsart}
\usepackage{pgf,tikz}
\usetikzlibrary{arrows}

\usepackage{undertilde}
\usepackage{mathtools}
\usepackage{amsthm}
\usepackage{amsmath}
\usepackage{verbatim}
\usepackage{amssymb}
\usepackage{graphicx}
\usepackage{latexsym}
\usepackage[mathscr]{euscript}

\newtheorem{theorem}{Theorem}[subsection]

\usepackage{color}

\newtheorem{lemma}[theorem]{Lemma}

\theoremstyle{definition}
\newtheorem{definition}[theorem]{Definition}

\theoremstyle{remark}
\newtheorem{remark}[theorem]{Remark}

\newcommand{\twon}{2^{\omega} }

\newcommand{\Lip}{\mathrm{\mathbf{Lip}} }

\newcommand{\zeroone}{[0,1] } 
\newcommand{\zn}{[0,1]^n }

\newcommand{\arpi}[2]{\Pi_{#2}^{#1}}
\newcommand{\arsigma}[2]{\Sigma_{#2}^{#1}}

\newcommand{\fcantor}{2^{<\omega}}

\newcommand{\dseq}[3]{ \left({#1}_{#2,#3}\right)_{{#2,#3}\in\NN} }

\newcommand{\real}{\mathbb{R}}

\newcommand{\lm}[1]{\lambda\left( #1 \right)}

\usepackage{ifpdf}

\ifpdf
\usepackage[bookmarks,colorlinks=true]{hyperref}
\else
\usepackage{hyperref}
\fi

\usepackage[UKenglish]{babel}

\newcommand{\id}{\text{id}}

\newcommand{\B}{\mathcal{B}}

\newcommand{\K}{\mathcal{K}}
\renewcommand{\L}{\mathcal{L}}

\newcommand{\N}{\mathbb{N}}

\renewcommand{\P}{\mathcal{P}}

\newcommand{\Q}{\mathbb{Q}}
\newcommand{\R}{\mathbb{R}}

\theoremstyle{plain}

\newtheorem{claim}[theorem]{Claim}

\newcommand{\NN}{{\mathbb{N}}}
\newcommand{\RR}{{\mathbb{R}}}

\newcommand{\uhr}[1]{\! \upharpoonright_{#1}}

\renewcommand{\P}{\mathcal P}
\newcommand{\bi}{\begin{itemize}}
\newcommand{\ei}{\end{itemize}}
\newcommand{\bc}{\begin{center}}
\newcommand{\ec}{\end{center}}

\newcommand{\strcantor}{2^{ < \omega}}

\renewcommand{\P}{\mathcal P}

\newcommand{\n}{\noindent}

\newcommand{\seq}[1]{(#1_i)_{i\in\NN}}

\newcommand{\ddos}[3]{\overline D_+{#1}\left({#2};{#3}\right)}
\newcommand{\ddoi}[3]{\underline D_+{#1}\left({#2};{#3}\right)}

\newcommand{\ddo}[3]{D_+{#1}\left({#2};{#3}\right)}
\newcommand{\dd}[3]{D{#1}\left({#2};{#3}\right)}

\newcommand{\vectornorm}[1]{\left|\left|#1\right|\right|}
\usepackage[normalem]{ulem}

\binoppenalty=\maxdimen
\relpenalty=\maxdimen

\begin{document}

\title{Randomness and differentiability in higher dimensions}

\author[\mbox{Galicki and Turetsky}]{Alex Galicki and Daniel Turetsky}

\address{A.~Galicki, Department of Computer Sciece, 
University of Auckland}
\email{agal629@aucklanduni.ac.nz}
\address{D. Turetsky, Kurt G\"odel Research Center}
\email{daniel.turetsky@univie.ac.at}
\begin{abstract}
We present two theorems concerned with algorithmic randomness and differentiability of functions of several variables. Firstly, we prove an effective form of the Rademacher's Theorem: we show that computable randomness implies differentiability of computable Lipschitz functions of several variables. Secondly,  we show that weak 2-randomness is equivalent to differentiability of computable a.e. differentiable functions of several variables.  
 \end{abstract}

\maketitle


\section{Introduction}

\subsection{Introduction}

The main subject of this paper lies at the interface of computable analysis (\cite{Weihrauch:00}) and algorithmic randomness (\cite{Nies:book}, \cite{Downey.Hirschfeldt:book}).

Intuitively, a real number is random if it does not have any exceptional properties. This approach can be formalized via identifying exceptional properties with effective null sets. To different types of effective null sets correspond different notions of  algorithmic randomness. 

One of the most fruitful areas of research concerned with interconnections between the two subjects is differentiability of effective functions. The main reason for this is that sufficiently well-behaved functions are almost everywhere differentiable. In this case the set of non-differentiability points of an effective function forms an effective null set and thus a test for algorithmic randomness. This makes it possible to characterize different randomness notions in terms of sets of points of differentiability of effective functions. Conversely, sets of points of differentiability for functions of particular classes can be characterised in terms of algorithmic randomness. The results of this kind are particularly compelling, since they show non-trivial connections between two seemingly distant areas of mathematics.

In recent years, a number of results of that kind have been published (for example, see \cite{Brattka.Miller.ea:nd,Miyabe:12, Pathak.Rojas.ea:12, Freer.Kjos.ea:nd}). Most of them are concerned with functions of one variable. Relatively few results are known about effective functions of several variables. 

Our first result is concerned with Lipschitz functions, which are particularly well behaved and enjoy a lot of attentions from mathematicians since they appear naturally in various contexts. The following classical result is called Rademacher's Theorem (see Section 3.1 in \cite{Evans.Gariepy:92}), it states that Lipschitz functions are almost everywhere\ differentiable.    

\begin{theorem}[Rademacher, \cite{Rademacher:19}] Suppose $U$ is an open subset of $\R^n$ and \mbox{$f:U\to\R^m$} is a Lipschitz function. Then there exists a null set, such that  $f$ is differentiable outside it.
\end{theorem}

\n We prove  the following effective form of Rademacher's Theorem. 

\newtheorem*{thm:rademacher}{Theorem \ref{t_rademacher_1}}
\begin{thm:rademacher}
Let $f:\zn\to\real$ be a computable Lipschitz function and let $z\in\zn$ be computably random. Then $f$ is differentiable at $z$.
\end{thm:rademacher}

The one dimensional variant of effective Rademacher's Theorem and its converse have been proven in \cite{Freer.Kjos.ea:nd}. 
\begin{theorem}[Theorem 4.2 in \cite{Freer.Kjos.ea:nd}]\label{thm:lipshitz1} A real $z\in\zeroone$ is computably random $\iff$ each computable Lipschitz function $f:\zeroone\to\R$ is differentiable at $z$.
\end{theorem}

Theorem \ref{t_rademacher_1} generalizes the $\Rightarrow$ direction of the above result.

The question whether the converse of the classical Rademacher's Theorem holds, that is whether every Lebesgue null-set is contained in a set of non-differentiability points  of a Lipschitz function, has been answered very recently after several decades of work by classical analysts (see \cite{Alberti.Csornyei.Preiss:10} and \cite{Preiss.Speight:14}). The converse holds when $m\ge n$ and does not hold otherwise.

To characterise differentiability sets of effective functions in terms of randomness in the usual way, those functions must be differentiable almost everywhere, for otherwise the sets of non-differentiability do not form null sets and cannot be interpreted as exceptional properties. This implies that the broadest possible class of functions in this context is the class of almost everywhere differentiable functions. For functions of one variable the following result is known:

\begin{theorem}[Theorem 6.1 in \cite{Brattka.Miller.ea:nd} ]\label{thm:w2r1}
Let $z\in \zeroone.$ The following are equivalent:
\begin{enumerate}
\item $z$ is weakly 2-random, and
\item all computable a.e.\ differentiable functions are differentiable at $z$.
\end{enumerate}
\end{theorem}  

\n Our final result is the the following generalization of that theorem.

\newtheorem*{thm:w2r}{Theorem \ref{w2r_theorem}}
\begin{thm:w2r}
Let $z\in\zeroone^n$, then the following are equivalent:
\begin{enumerate}
\item $z\text{ is weakly $2$-random}$,
\item all partial derivatives exist for all computable a.e.\ differentiable \mbox{$f:\zn\to\real$},
\item each computable a.e.\ differentiable function is differentiable at $z$.
\end{enumerate}
\end{thm:w2r}

\subsection{Structure of the paper}
\n In the rest of this section we present  relevant definitions and facts  and introduce some useful notation. 

In Section 2 we prove an effective version of Rademacher's Theorem. We start the section by recalling some important facts about Lipschitz functions and then proceed with the proof of the main result. The section ends with a discussion of a relatively recent classical result of Maleva and Dor\'e \cite{Dore.Maleva:11} and some of its implications.

In Section 3 we demonstrate that weak 2-randomness characterises differentiability points of computable a.e.\ differentiable functions. 

The last section discusses some open problems related to this article.  

\subsection{Preliminaries}

\subsubsection{Measure}

We work exclusively with the Lebesgue measure on $\zn$. Slightly abusing notation, we always denote it by $\lambda$.

\subsubsection{Derivatives in higher dimensions} Let $f:\zn\to\real$ be a function and let $x\in\zn$. We say $f$ is \emph{differentiable} at $x$ if for some linear map $T$ the following holds $$\lim_{h\to 0}\frac{f(x+h)-f(x)-T\cdot h}{\vectornorm h}=0.$$ Then, by definition, $f'(x)=T$.

\n Let $\{e_i:1\le i\le n\}$ denote the standard basis for $\RR^n$. We denote \emph{partial derivatives} by $D_if(x)$, lower and upper partial derivatives by $\underline D_if(x)$ and $\overline D_if(x)$, respectively.

Working with derivatives often means working with slopes. In our proofs we use the following notation.

\n Fix coordinate~$i$.  For $x \in \zn$, $h \in \RR$, define
\[
\delta_f^i(x, h) = \frac{f(x_1, \dots, x_i + h, \dots, x_n) - f(x_1, \dots, x_i, \dots, x_n)}h,
\]

\n
and 
\[
\delta^{1..n}_f({x}, h) = \left[\begin{array}{ccc} \delta_f^1({x},h) & \dots & \delta_f^n({x},h) \end{array}\right].
\]

\subsubsection{Computable real functions} There are multiple ways of formalizing computability of real functions, most of which turned out to be equivalent. We will rely on the following definitions.

A sequence $\seq q$ of elements of $\R^n$ is called a \emph{Cauchy name} if the coordinates of each $q_i$ are rational, and $\|q_k-q_n\|\le 2^{-n}$ for all $n,k$ with $k\ge n$.
If $\lim_{n\to\infty} q_n=x$, then we say that $\seq q$ is a \emph{Cauchy name for} $x$.

We say $x\in \R^n$ is \emph{computable} if there is a computable Cauchy name for $x$.

\begin{definition}
A function $f:\zn\to\real^m$ is \emph{computable} if:
\begin{enumerate}
\item $f(q)$ is computable (uniformly in $q$) where $q$ has all dyadic rational coordinates, and

\item $f$ is \emph{effectively uniformly computable}, that is if there is a computable \mbox{$h:\N\to\N$} such that $\|x-y\|\le 2^{-h(i)}$ implies $\|f(x)-f(y)\|\le 2^{-i}$ for all $x,y\in\zn$ and all $i\in\N$.
\end{enumerate}
\end{definition}

A more intuitive understanding of the above definition is that $f$\ is computable if there is an algorithm that, given a Cauchy name for $x$, computes a Cauchy name for $f(x)$.

\subsubsection{Algorithmic randomness}~

\n The most common method for defining a randomness notion is via effective null sets.  The following two randomness notions are of direct interest to us and are defined in terms of avoidance of effective null sets.

\begin{definition} 
Let $z\in\zn$. We say $z$ is \emph{weakly random} if there does not exist a $\arpi01$ null set that contains $z$. 
Similarly, we say $z$ is \emph{weakly 2-random} if there does not exist a $\arpi02$ null set that contains $z$.
\end{definition} 

An alternative approach to formalizing randomness notions is via effective betting strategies. An infinite binary string can be thought of as random if no (effective) betting strategy can succeed by betting on bits of that string. Betting strategies are usually formalized as \emph{martingales} (see \cite{Nies:book}, Chapter 7).

\begin{definition}
We say a function $B:\fcantor\to\Q_+$ is a \emph{martingale} if the following condition holds for all $\sigma\in\fcantor:$ $$2B(\sigma)=B(\sigma0)+B(\sigma1).$$

\n $B(\sigma)$ can be interpreted as the value of capital after betting on bits of $\sigma$.
We say $B$ \emph{succeeds} on $Z\in\twon$ if $\liminf_n B(Z\uhr n)=\infty$.
\end{definition}

\begin{definition}\label{cr_definition}
We say $Z\in\twon$ is \emph{computably random} if no computable martingale succeeds on $Z$.

We say $z=(0.Z_1,\dots,0.Z_n)\in\zn$ is computably random if its binary expansion, that is $Z=Z_1\oplus\dots\oplus Z_n$, is computably random. Here $0.A$ denotes the real number whose binary expansion is $A\in\twon$.
\end{definition}

\n It is known that weak 2-randomness implies computable randomness and computable randomness implies weak randomness.

\subsubsection{Preservation of computable randomness}

\begin{definition}(cf. 7.1 in \cite{Rute:13})

We say that $\phi:\zn\to\zn$ is \emph{almost everywhere (a.e.)\ computable} if there exists a partial computable $F:\N^\omega\to\N^\omega$ and a $\arpi02$ subset $A\subseteq\zn$ with $\lm A=1$ such that:
\begin{enumerate}
\item  for all $x\in A$, given a Cauchy name of $x$, $F$ computes a Cauchy name for $\phi(x)$, and
\item $x\in A$ iff for all $a,b$, which are Cauchy names for $x$, both $F(a)$ and $F(b)$ are Cauchy names for the same element.
\end{enumerate}

We say that $\phi:\zn\to\zn$ is an \emph{a.e.\ computable isomorphism} if there exists $\psi:\zn\to\zn$ such that $\phi\circ\psi=\id$ and $\psi\circ\phi=\id$ almost everywhere and both $\psi,\phi$ are measure preserving and a.e.\ computable.
\end{definition}

We are interested in the above notions for the following property of computable randomness proven by Rute.

\begin{theorem}[Theorem 7.9 in \cite{Rute:13}]\label{cr_preservation}
Let $T$ be an a.e.\ computable isomorphism. Then    for all $x\in\zn$, $x$ is computably random if and only if $T(x)$ is computably random. 
\end{theorem}

\subsubsection{Uniform relative computable randomness}~

\n\ Both the following definition and theorem are  due to Miyabe and Rute (\cite{Miyabe.Rute:13}).

\begin{definition}
A total computable function $m:\twon\times\strcantor\to \real$ is a \emph{uniform computable martingale} if $m(Z,\cdot)$ is a martingale for every $Z\in\twon$.

We say $A$ is \emph{computably random uniformly relative to }$B$ if there is no uniform computable martingale $m$ such that $m(B,\cdot)$ succeeds on $A$.

\end{definition}

\n Note that the above definition works for elements of $\zn$ as well. 
\begin{theorem}[Theorem 1.3 in \cite{Miyabe.Rute:13}]\label{t_miyabe_rute1}
$A\oplus B$ is computably random if and only if $A$ is computably random uniformly relative to $B$ and $B$ is computably random uniformly relative to $A$.
\end{theorem}

\section{Effective form of Rademacher's Theorem}\label{lip_section}

\n In this section we prove a theorem which can be seen as an effective version of Rademacher's.

\begin{theorem}\label{t_rademacher_1}
Let $f:\zn\to\real$ be a computable Lipschitz function and let $z\in\zn$ be computably random. Then $f$ is differentiable at $z$.
\end{theorem}

\begin{remark}
An immediate consequence of the above theorem is that computable randomness of $z\in\zn$ is sufficient for differentiability of computable Lipschitz functions form $\zn$ to $\R^m$ for any $n,m$.
\end{remark}

\n Lipschitz functions are particularly well-behaved and have a number of properties related to differentiability in general and to directional derivatives in particular. Some of those properties will be used by us in the proof of the above theorem and this is why we start this section by recalling some facts about Lipschitz functions and by establishing some useful notation before proceeding to the proof.

\subsection{Lipschitz functions}~

A function $f:\R^n\to\R^m$ is \emph{Lipschitz} if there exists $L\in\R^+$ such that
\[\|f(x)-f(y) \|\le L\cdot\|x-y\|\text{ for all }x,y\in\R^n.\]

\n The least such $L$ is called \emph{the Lipschitz constant} for $f$. We denote it by $\Lip(f)$.

\vspace{5pt}

 Let $K_n\subset\real^n$ be defined as $K_n=\{(x_1,\dots,x_n)\in\real^n:x_i\ge 0 \text{ for all } i\le n\}$. We say $f:\R^n\to \R$  is \emph{$K_n$-increasing} if $f(x+k)\ge f(x)$ for all $k\in K_n$. $f$ is called \emph{$K_n$-monotone} if either $f$ or $-f$ is $K_n$-increasing. 

\begin{remark}\label{decomposition_remark}
\n Every Lipschitz function $f:\R^n\to\R$ is a sum of two $K_n$-monotone functions. To see this, let $ {m}=(\Lip(f),\dots,\Lip(f))\in\R^n$ and note that\\ $f=(f+\langle m, x\rangle)-\langle m, x\rangle$, and that both summands are $K_n$-monotone.
\end{remark}

\n 

\subsection{Directional, G\^ateaux and Fr\'echet derivatives}~

\n In order to exploit some of the properties of Lipschitz functions, we need to present a more nuanced view of differentiability in higher dimensions.

\n Let $f:\real^n\to\real$ be a function, we define the \emph{Dini-directional deri\-vatives}
of $f$ at a point $x\in \R^n$ with respect to a direction $v\in\R^n$ as 
\begin{align*}
\ddos fxv=\limsup_{t\downarrow 0}\frac{f(x+tv)-f(x)}{t}\text{ and }\ddoi fxv=\liminf_{t\downarrow 0}\frac{f(x+tv)-f(x)}{t}.
\end{align*}

\n When $\ddos fxv=\ddoi fxv$ is finite, we define \emph{one-sided directional derivative} by $$\ddo fxv=\lim_{t\downarrow 0}\frac{f(x+tv)-f(x)}{t}.$$

\n The \emph{two-sided directional derivative} $\dd fxv$ is defined by $$\dd fxv=\lim_{t\to 0}\frac{f(x+tv)-f(x)}{t}.$$

\n To work with directional slopes, we need the following notation.

\n For $x \in \zn$, $v\in\R^n$, and  $h \in \RR$, define
\[
\delta_f^v(x, h) = \frac{f(x+hv) - f(x)}h.
\]

\vspace{10pt}
\n If all two-sided directional derivatives of $f$ at $x$ exist and the function $T$ given by $T(y)=\dd fxy$ is linear, then $f$ is said to be \emph{G\^ateaux-differentiable} at $x$. The linear map $T$ is called the \emph{G\^ateaux derivative} of $f$ at $x$. Furthermore, if $f$ is G\^ateaux-differentiable at $x$ and if $$\lim_{h\to 0}\frac{f(x+h)-f(x)-T\cdot h}{\vectornorm h}=0,$$ then $f$ is said to be \emph{Fr\'echet differentiable} at $x$.

\n Thus, Fr\'echet differentiability is equivalent to the usual differentiability.

The following observation is crucial
for the main proof and justifies presenting differentiability in this more elaborate way.

\begin{remark}\label{lip_remark}

\n  
For Lipschitz functions on $\R^n$, G\^ateaux and Fr\'echet differentiability coincide (for example, see Observation 9.2.2 in \cite{Lindenstrauss.Preissr:12}).

Furthermore, it is known (see \cite{Chabrillac.Crouzeix:87}) that for a $K_n$-monotone function $f$, both $\ddos fx\cdot$ and $\ddoi fx\cdot$ are continuous on the interior of $K_n\cup ~{-K_n}$. Thus, for a Lipschitz function $f$, both $\ddos fx\cdot$ and $\ddoi fx\cdot$ are continuous everywhere.

The following property is  a direct consequence of the above fact. Let $A$ be a dense subset of $\real^n$, let $f:\zn\to\real$ be a Lipschitz function and let $x\in\zn$, then: 
\begin{enumerate}
\item[$(\star)$] if $v\mapsto \ddo fxv$ is defined and is linear on $A$, then $v\mapsto \ddo fxv$ is defined everywhere and is linear.

\end{enumerate}

\n This means that in order to show that a Lipschitz function $f$ is differentiable at~$x$, it is sufficient to show that  $\dd fx\cdot$ is defined and linear on a dense subset of directions. 

\end{remark}

\subsection{Overview of the proof.}~ 

\n The proof consists of three distinct steps.
 
(1) We show that all partial derivatives of $f$ at $z$ exist. Firstly, we prove an analogous result for $K_n-$monotone functions and then use the fact that every Lipschitz function is a sum of two $K_n-$monotone functions to prove the required result holds for Lipschitz functions. The result for $K_n-$monotone functions is a consequence of the following two facts: (i) a uniform relativization of the $\Rightarrow$ implication of Theorem \ref{t_nies1} and (ii) a form of Van Lambalgen's Theorem for computable randomness proven by Miyabe and Rute \cite{Miyabe.Rute:13}.   

(2) We use the above fact to show that 
 all one-sided directional derivatives of $f$ at $z$ exist. Since, by Remark \ref{lip_remark}, we are  only required to show this for a dense set of directions, we only consider computable directions $v$. Two observations play a crucial role in this step: (a) computable randomness is preserved by computable linear isometries, and (b) linear functions are Lipschitz. The above observations are used to define a computable Lipschitz function $g_v:\zn\to\R$ such that $D_1g_v(\hat z)=\dd fzv$ for some computably random $\hat z\in \zn$. By the result proven in the first step, $D_1g_v(\hat z)$ does exist.
 
(3) Finally, we show that the function  $T(u)=\ddo fzu$ is linear. Again, we consider only computable directions. We show that any point where directional derivative is not linear and the failure of linearity is witnessed by a computable direction, belongs to a $\arpi01$ null set. Since $z$ is computably random, this completes the proof.

\n Showing the linearity of $\ddo fz\cdot$ is the final step of our proof, because for Lipschitz functions, G\^ ateaux differentiability implies (full)\ differentiability.

\subsection{Existence of partial derivatives} Firstly, we show that computable randomness is sufficient for all partial derivatives of computable $K_n-$monotone functions to exist. An analogous result for computable Lipschitz functions is a simple corollary of that.

To achieve the required result, we combine a variation of Van Lambalgen's Theorem for computable randomness, proven by Miyabe and Rute (\cite{Miyabe.Rute:13}),  with a variation of the $\Rightarrow$ implication of the following result.

\begin{theorem}[Theorem 4.1 in \cite{Brattka.Miller.ea:nd}]\label{t_nies1} A real $x$ is computably random $\iff$ $f'(x)$ exists for each computable nondecreasing function $f:\zeroone\to\real$.
\end{theorem}

\begin{lemma}\label{uniform_relativization_monotone} 

Let $g:\twon\times \zeroone\to \R$ be a total computable function such that $g(X,\cdot)$ is monotone for all $X\in\twon$ and let $Z,Y\in\twon$. If $Z\oplus Y$ is computably random, then $g_Y'(z)$ exists, where $g_Y=g(Y,\cdot)$ and $z=0.Z$.
\end{lemma}

The proof of Lemma \ref{uniform_relativization_monotone} is a modification of the proof of Theorem 4 in \cite{nies:2014}. Theorem 4 in \cite{nies:2014} is a polynomial version of Theorem \ref{t_nies1}, but its proof is somewhat simpler and requires only a few modifications to yield the kind of uniform relativization needed for the proof of Lemma \ref{uniform_relativization_monotone}. We will describe the required changes without repeating the whole proof.

\begin{proof}

The proof is by contraposition. Let $g:\twon\times \zeroone\to \R$ be a total computable function such that $g(X,\cdot)$ is monotone for all $X\in\twon$. Let $z\in\zeroone$, let $Z$ be the binary expansion of $z$ and let $Y\in\twon$.
Define $g_Y=g(Y,\cdot)$ and suppose $g_Y'(z)$ doesn't exist.
We need to exhibit a uniform computable martingale $d$ such that $d(Y,\cdot)$ succeeds on $Z$. 

In the $\Rightarrow$ direction of the original proof, assuming $f'(x)$ does not exists (where $f:\zeroone\to\R$ is a polynomial time computable monotone function), Nies constructed a (polynomial time) computable martingale that succeeds on the binary expansion of $x$. The assumption that $f$ is polynomial time computable was used to show that the resulting martingale is polynomial time computable. If this assumption is relaxed so that $f$ is assumed to be computable, the resulting martingale ends up being computable, rather than polynomial time computable. 
We need to verify that a slightly modified proof works for demonstrating that there is  a uniform computable martingale $d$ such that $d(Y,\cdot)$ succeeds on $Z$.

\emph{Uniform relativization of the $\Rightarrow$ implication of Theorem 4 in \cite{nies:2014}.} Here we use the combined terminology from the original proof and the terminology required for our proof. For the $\Rightarrow$ direction of the proof of Theorem 4 in \cite{nies:2014}, Nies had to consider two cases: $\widetilde D_2f(x)<\widetilde D f(x)$ and $\utilde Df(x)<\utilde D_2 f(x)$. Nies constructed a pair of computable martingales, $L$ and $L'$, corresponding to the above mentioned cases, such that either $L$ succeeds on the binary expansion of $x$, or $L'$ succeeds on the binary expansion of $x-1/3.$ Both $L$ and $L'$ query the same martingale $M$ defined by $M(\sigma)=S_f([\sigma])$.
Since $M:\twon\times\fcantor\to\R$ defined by $$M(Y,\sigma)=S_{g(Y,\cdot)}([\sigma])$$ is a uniform computable martingale, it can be easily checked that constructions of $L$ and $L'$ can be naturally extended to define uniform computable martingales $\L$ and $\L'$ such that either 
\begin{enumerate}
\item $\L(Y,\cdot)$ succeeds on $Z$ or 
\item $\L'(Y,\cdot)$ succeeds on the binary expansion of  $z-1/3$ (without loss of gene\-rality we may assume that $z>1/3$).
\end{enumerate}
The first case implies that $Z$ is not computably random uniformly relative to $Y$ and thus $Z\oplus Y$ is not computably random.
 
Note that $(x_1,x_2)\mapsto (x_1,x_2+1/3\mod 1)$ is an a.e.\ computable isomorphism. And since computable randomness is preserved by a.e.\ computable isomorphisms, the second case  implies that $Z\oplus Y$ is not computably random.  

\end{proof}

\begin{remark}The original proof relied on a different preservation property of computable randomness. It was using the fact that computable randomness is base invariant. We could not use the result about base invariance in our proof immediately (since we have now multiple coordinates instead of one), hence we chose to use another preservation property of computable randomness.
\end{remark}

\n

\begin{lemma}\label{l_cr_implies_partial_derivatives}

Let $z\in\zn$ be computably random and let $f:\zn\to\real$ be a computable $K_n-$increasing function. Then all partial derivatives of $f$ at $z$ exist.

\begin{proof}

Fix $i\le n$. The proof is by contraposition: suppose $D_if(z)$ does not exist, we will show that $z$ is not computably random.

Let $y=z-z_ie_i$ and let $Y$ be its binary expansion. Define $g:\twon\times\zeroone\to\R$ by  $$g(X,h)=f(0.X+he_i)$$ and let $g_y=g(Y,\cdot)$. Then $g$ satisfies relevant assumptions of Lemma \ref{uniform_relativization_monotone} and  $g_y'(z_i)=D_if(z)$. Furthermore, we know that $g_y'(z_i)$ does not exist. To show $z$ is not computably random, by Theorem \ref{t_miyabe_rute1}, it is sufficient to  show that $Z_i$ is not computably random uniformly relative to $Y$ (as this implies $Z_i$ not being computably random uniformly relative to $\oplus_{j\neq i} Z_j$). This follows from Lemma \ref{uniform_relativization_monotone}.
\end{proof}
\end{lemma}

\begin{lemma}\label{l_cr_implies_partial_derivatives_lipschitz}

Let $z\in\zn$ be computably random and let $f:\zn\to\real$ be a computable Lipschitz function. Then all partial derivatives of $f$ at $z$ exist.

\begin{proof}
Similar to the Remark \ref{decomposition_remark}, let $M=\text{Lip}(f)$ and let $\bold m=(M,\dots,M)\in\real^n$, then $g(x)=f(x)+\langle \bold m, x\rangle$ is a $K_n-$increasing computable function. Thus all partial derivatives of $g$ at $z$ exist, and therefore all partial derivatives of $f$ at $z$ exist too.
\end{proof}
\end{lemma}

\subsection{Existence of directional derivatives}\label{existence_directional_derivatives}

We will use the previously proven fact about existence of partial derivatives of Lipschitz functions  to show that, in fact, an analogous result holds for all one-sided directional derivatives. The main idea relies on two simple observations: 
\begin{enumerate}
\item computable randomness is invariant under computable linear isometries, and

\item linear functions are Lipschitz. 

\end{enumerate}

\vspace{5pt}
\n For any $u,v\in\real^n$ with $\|v\|=\|u\|=1$ and $u\neq v$, fix (say, via the Gram-Schmidt process) two orthonormal bases $B_u,B_v$ of $\R^n$ with $v\in B_v$ and $u\in B_u$.  Let \mbox{$\Theta_{u\to v}:\real^n\to\real^n$} denote a change of basis map (that takes $B_u$ to $B_v$) such that $\Theta_{u\to v}(u)=v$. This function is a linear isometry and it is computable when $u,v$ are computable. 

\n The image of the unit cube $\zn$ under functions of the form \mbox{$\Theta_{u\to v}:\real^n\to\real^n$} is not necessarily contained in $\zn$. To deal with this issue, we use the function $\P_1:\real^n\to\zn$ defined by  $$\P_1(x_1,\dots,x_n)=(\min\{1,x_1\},\dots, \min\{1,x_n\}).$$
$\P_1$ is a computable Lipschitz function which coincides with the identity map on the unit $n$-cube.
\n For any function $f:\zn\to\real$, let $\hat f=f\circ \P_1$, so that if $f$ is computable and a.e.\ differentiable, so is $\hat f$. Moreover, if $f$ is Lipschitz, so is $\hat f$. Note that $\hat f$ is defined on the whole $\real^n$.

\begin{lemma}\label{lemma_directional_derivative_trick}
\n Let $f:\real^n\to \real$ be a function, let $u,v,w\in\real^n ,x\in\zn$ and let $\Theta=\Theta_{v\to u}$.  Then $$\ddo fxu = \ddo gzv$$ where $g=f\circ (\Theta+w)$ and $z=\Theta^{-1}(x-w)$.

\begin{proof}~

\n First, note that for any $t>0,$
\begin{align*}
\frac{g(z+tv)-g(z)}{t}=\frac{f(\Theta(z+tv)+w)-f(\Theta(z)+w)}{t}=\frac{f(x+tu)-f(x)}{t}.
\end{align*}

\n By taking the limits of both sides we get the required equality. 

\end{proof}
\end{lemma}

\begin{lemma}\label{lemma222}
Let $f:\zn\to\real$ be computable Lipschitz and suppose $x\in\zn$ is computably random. Then $\ddo fxu$ exists for every $u\in\real^n$.

\begin{proof}
It is sufficient to show that $\ddo fxu$ exists for each computable $u$ with $\|u\|=1$.  

Let $u$ be computable and let $v=e_1$. By density we can find some computable $w\in\real^n$, so that $z=\Theta_{v\to u}^{-1}(x-w)$ is contained in $\zn$.

We apply Lemma \ref{lemma_directional_derivative_trick} to $\hat f,v,u,w$ and $x$, so that $$\ddo fxu=\ddo{\hat f}xu=\ddo gzv$$ where $g$ is Lipschitz and computable and $z\in\zn$ is computably random (again, we use Theorem \ref{cr_preservation} here). The required result follows from the fact that $\ddo gzv=D_1g(z)$ and we know that $D_1g(z)$ exists.

\end{proof}
\end{lemma}

\subsection{Linearity of directional derivatives}~

\n In the last step of the proof, we need to show that $\ddo fz\cdot$ is linear on computable elements (where $f$ is computable Lipschitz and $z$ is computably random).

\n Let $f:\zn\to\real$ be a function. For $u\in\real^n$, define

$$\K^f_u=\{z~|~ \ddo fzu \text{ exists} \}.$$ 

\n For $q\in\Q^{+}$ and $u,v\in \real^n$, define $\L_{u,v,q}^f$ to be 
the set of points where linearity of $\dd fz\cdot$ fails and the failure is witnessed by $u,v$ and $q$. More formally, let 
 $$\L_{u,v,q}^f=\K^f_u\cap \K^f_v\cap \K^f_{u+v}\cap \{z~|~ |\ddo fz{u+v}- \ddo fzu- \ddo fzv|\ge q\}.$$ 

\begin{lemma}\label{lemma302} Let $f:\zn\to\real$ be a computable a.e.\ differentiable function. Let $v,u\in\real^n$ be computable.
Let $z\in\L^f_{v,u,q}$ for some $q\in\Q$. Then there exist a $\arpi01$ null-set that contains $z$. 

\begin{proof} 

Since $\ddo fzv,\ddo fzu$ and $\ddo fz{v+u}$ exist, there is $p>0$ such that 
$\left| \delta^v_f(z,h)+\delta^u_f(z,h)-\delta^{v+u}_f(z,h)  \right|\ge q$ for all $h\le p$. Hence the set of all $x$ such that $$\forall h\, \left(h  \le p \implies \left| \delta^v_f(x,h)+\delta^u_f(x,h)-\delta^{v+u}_f(x,h)  \right|\ge q \right),$$ where $h$ range over rationals, contains $z$. It is clearly a $\arpi01$ set and it is a null set, since its complement contains all  points of differentiability of $f$ and  $f$ is a.e.\ differentiable. 
\end{proof}

\end{lemma}

\n So far, we have shown that computable randomness implies existence of directional derivatives in and weak randomness is sufficient for linearity of directional derivatives. This implies that computable randomness is sufficient for G\^ateaux differentiability and this completes the proof of Theorem \ref{t_rademacher_1} since G\^ateaux differentiability implies differentiability of Lipschitz functions on $\zn$.

\subsection{Compact universal null-sets}\label{dore_maleva}

In the context of differentiability of Lipschitz functions, a subset $A$ of $\R^n$ is said to be \emph{universal} if every real-valued Lipschitz function on $\R^n$ is differentiable at some point of $A$. Since the early work od {D.~Preiss~\cite{Preiss:90},} it is known that there exist universal $G_\delta$ null sets for $n\ge 2$. Relatively recently, Dor\'e and Maleva constructed a family of compact universal null sets (see \cite{Dore.Maleva:11}). The crucial idea in their construction is that a Lipschitz function is differentiable at points where a directional derivative  is maximal in some specific sense and that such points can be found on small line segments (see Lemmata 4.2 and 4.3 in \cite{Dore.Maleva:11}). Their sets contain lots of such line segments and this is the reason they are universal. The result implies (as will be shown shortly) the existence of a universal $\arpi01$ null set and this has significant implications for tackling the question of which randomness notion is implied by the Rademacher's Theorem. 

To characterise a randomness notion $X\subset\R^n$ via differentiability of computable real-valued Lipschitz functions, it is sufficient to prove two statements: 
\begin{enumerate}
\item $z\in X$ implies differentiability of all computable real-valued Lipschitz functions at $z,$ and 
\item 
differentiability of all computable real valued Lipschitz functions at $z$, implies $z\in X$. 
\end{enumerate}
The second type of statements is usually proven by explicitly constructing a computable function not differentiable at a given randomness test (for $X$). The existence of a universal $\arpi01$ null set shows that such an approach cannot succeed even in proving that differentiability of computable real-valued Lipschitz functions implies weak randomness. 

The construction in \cite{Dore.Maleva:11} is parameterized by two sequences. 
Below we verify that with suitable parameters this construction yields a $\arpi01$ null set.

\n\emph{Construction by Dor\'e and Maleva.} 

Let $\seq N$ be a sequence of odd integers such that $N_1>1$, $\lim _i N_i= \infty$
 and $\sum \frac{1}{N^2_i}=\infty$. Let $\seq p$ be a sequence of real numbers with $1\le p_i\le  N_i$ and $\lim_i p_i/N_i=0$. Let $d_0=1$ and for all $i\ge 1$ let $d_i=\prod_{k\le i} N_k^{-1}$ and define a lattice in $\R^2$ $$C_i=\left(\frac{d_{i-1}}2,\frac{d_{i-1}}2\right)+\mathbb{Z}^2.$$

\n Finally, define $$W=\R^2\setminus\bigcup_{i\ge 1}\bigcup_{c\in C_i} \B_\infty(c, p_id_i/2),$$ where $\B_\infty(x,r)$ denotes an open ball in $\left(\R^2,\|\cdot\|_\infty\right)$.

$W$ is a closed null set. Dor\'e and Maleva proved [Corollary 5.2 in \cite{Dore.Maleva:11}] that for any such $W$, any open neighbourhood of  the set \mbox{$M=\R^{n-2}\times W$} contains a point of differentiability of every Lipschitz function $f:\R^n\to\real$. In particular, $\zn\cap M$ contains a point of differentiability of every Lipschitz $f:\zn\to\real$.

It is easy to see that both $\seq N$ and $\seq p$ can be taken to be computable sequences and then $\zn\cap M$ is a $\arpi01$ null set. (For example, take $\seq N$ to be $3,3,3,5,5,5,5,5,7,7,7,7,7,7,7,\dots$ and let $p_i=4$ for all $i$).

\section{Characterizing weak 2-randomness in terms of differentiability}

This section is devoted to proving Theorem \ref{w2r_theorem}, which characterises weak randomness in terms of differentiability of computable functions of several variables. It is worth pointing out that while our result is a generalization of  Theorem \ref{thm:w2r1}, it is ``stronger'' in the sense that we show equivalence of three conditions, rather than two. 
Recall that Theorem \ref{thm:w2r1} shows weak 2-randomness is equivalent to differentiability of computable a.e. differentiable functions. Somewhat surprisingly, in higher dimensions, a seemingly weaker condition, the existence of all partial derivatives for all computable a.e. differentiable functions, is also equivalent to weak 2-randomness.   

We start the section with a fact of independent interest.

\begin{lemma}\label{pi03_lemma} Let $f:\zn\to\real$ be a computable function. The set of points  at which $f$ is differentiable is a $\arpi03$ set.
\begin{proof}

Recall that the definition of the derivative of a function of several variables involves a nested limit. The main idea of this proof is that the set of points of differentiability (for a given function), $D$, can be written as an intersection of two sets, each of which can be described with only one limit.
Specifically, we will show that $D$  is the intersection of two $\arpi03$ sets, $A$ and $B$, where $A$\ is the set of points where all partial derivatives of $f$ exist, and $B$ is the set consisting of those~$x$ satisfying

\begin{align}\label{B_condition}
\lim_{||{h}|| \to 0,b\to 0} \frac{f({x} + {h}) - f({x}) - \delta^{1..n}_f({x}, b)\cdot{h}}{||{h}||} = 0.
\end{align}

\begin{claim}  $A$ is a $\arpi03$ set that contains all points of differentiability of $f$.
\begin{proof}
Fix coordinate~$i$. 
\n 
For~$q$ a rational, $\ \overline{D}_if(x) \geq q$ is equivalent to the formula
\[
\forall p  \, \forall \delta \, \exists h~ \left( |h|  < \delta \land \left(  (p<q\land \delta>0)\implies\delta_f^i(x,h) > p\right)\right).
\]
By density, we can take~$p$ and~$\delta$ to range over the rationals.  By continuity of~$f$, we can take~$h$ to range over the rationals.  Thus $\{ {x} \mid \overline{D}_if({x}) \geq q\}$ is a $\Pi^0_2$ set  uniformly in~$q$.  Symmetrically, so is $\{x \mid \underline{D}_if(x) \leq q\}$.  Then the set of $x$ such that $D_if(x)$ does not exist,  is precisely the set\[
\left\{ x : \forall q[\overline{D}_if(x) \geq q] \vee \forall q[\underline{D}_if(x) \leq q] \vee \exists q \exists p[\underline{D}_if(x) \leq q < p \leq \overline{D}_if({x})]\right\}.
\]
This is a $\Sigma^0_3$ set and hence the set of points where at least on partial derivative does not exist is also $\arsigma03$.  Thus $A$ is a $\arpi03$ set. The other part of the claim is trivial.

\end{proof}
\end{claim}

\

\begin{claim} $B$ is a $\arpi03$ set that contains all points of differentiability of $f$.

\begin{proof} 

By definition, $\lim_{|h| \to 0} \delta^{1..n}_f({x},h) = J_f({x})$, the Jacobian of~$f$ at~${x}$ (when this exists).

\n The derivative of~$f$ exists at~${x}$ if $J_f({x})$ exists and
\[
\lim_{||{h}|| \to 0} \frac{f({x} + {h}) - f({x}) - J_f({x})\cdot{h}}{||{h}||} = 0.
\]

\n\ To see that $B$ is a $\arpi03$ set, we can rewrite the condition (\ref{B_condition}) in the following form:

\[
\forall  \epsilon  \exists \delta \forall h \forall b  ~\left( \delta>0\right)\land\left[  \left( \epsilon > 0  \land ||{h}|| < \delta\land |b| < \delta\right)   \implies  \frac{\left| f(x + {h}) - f({x}) - \delta^{1..n}_f({x}, b)\cdot{h}\right| }{||{h}||} \le \epsilon\right].
\]
Here $\epsilon, \delta,h$ and~$b$ are rationals, and ${h}$ has rational coordinates.

Suppose~${x}$ is a point at which~$f$ is differentiable.  Fix $\epsilon > 0$.  Let $\delta$ be sufficiently small that
for all $h$ with $||{h}|| < \delta$, \[\frac{\left|f({x} + {h}) - f(x) - J_f({x})\cdot{h}\right|}{||{h}||} < \epsilon/2,
\]
and also
for all $b$ with  $|b| < \delta$,  

\[
\vectornorm{ (\delta^{1..n}_f({x},b) - J_f({x}))^T } < \epsilon/2.
\]
Here we treat $\delta^{1..n}_f({x},b) - J_f({x})$ as a row vector.  Then for any $h$ and $b$ with $||{h}|| < \delta$ and $|b| < \delta$,
\begin{eqnarray*}
\frac{ \left| f(x + {h}) - f(x) - \delta^{1..n}_f(x, b)\cdot{h}\right| }{||{h}||}
&=& \frac{\left|f({x} + {h}) - f({x}) - J_f({x}){h} + J_f({x}){h} - \delta^{1..n}_f({x}, b)\cdot{h}\right|}{||{h}||}\\
&\leq& \frac{\left|f({x} + {h}) - f({x}) - J_f({x})\cdot{h}\right|}{||{h}||}
        + \frac{\left|(J_f({x}) - \delta^{1..n}_f({x},b))\cdot{h}\right|}{||{h}||}\\
&<& \epsilon/2 + \frac{\vectornorm{(J_f({x}) - \delta^{1..n}_f({x},b))^T}||{h}||}{||{h}||}\\
&<& \epsilon/2 + \epsilon/2 = \epsilon.
\end{eqnarray*}
Thus~$B$ contains every point at which~$f$ is differentiable.

\end{proof}
\end{claim}

\n Thus, $A\cap B$ contains all points of differentiability of $f$. Let's show that the converse inclusion holds.\ 

\begin{claim}
$f$ is differentiable at all elements of $A\cap B$.
\begin{proof}

Let $x\in A\cap B$. Fix~$\epsilon > 0$.  Since $x\in B$, we can find~$\delta$ such  that
\[
\forall {h}\forall b\left[\left(|b| < \delta \land ||{h}|| < \delta \right) \implies \frac{\left|f(x + {h}) - f(x) - \delta^{1..n}_f(x,b)\cdot{h}\right|}{||{h}||} < \epsilon/2\right],
\]
and since all partial derivatives of $f$ at $x$ exist, we can find some $b$ with $|b| < \delta$ such that
\[
\vectornorm{ (\delta^{1..n}_f(x,b) - J_f(x))^T } < \epsilon/2.
\]
Then, for any $h$ with $||{h}|| < \delta$,
\begin{eqnarray*}
\frac{ \left| f(x + {h}) - f(x) - J_f(x)\cdot{h}\right|}{||{h}||}
&=& \frac{ \left| f(x + {h}) - f(x) - \delta^{1..n}_f(x,b)\cdot{h} + \delta^{1..n}_f(x,b)\cdot{h} - J_f(x)\cdot{h}\right|}{||{h}||}\\
&\leq& \frac{ \left| f(x + {h}) - f(x) - \delta^{1..n}_f(x,b)\cdot{h}\right|}{||{h}||} + \frac{ \left| (\delta^{1..n}_f(x,b) - J_f(x))\cdot{h} \right|}{||{h}||}\\
&<& \epsilon/2 + \epsilon/2 = \epsilon.
\end{eqnarray*}
Thus~$f$ is differentiable at~$x$.
\end{proof}
\end{claim}

\end{proof}
\end{lemma}

\begin{theorem}\label{w2r_theorem}
Let $z\in\zeroone^n$. The following are equivalent:
\begin{enumerate}
\item $z\text{ is weakly $2$-random,}$
\item all $D_if(z)$ exist for all computable a.e.\ differentiable $f:\zn\to\real,$
\item each computable a.e.\ differentiable function is differentiable at $z.$
\end{enumerate}
\begin{proof}[Proof (1) $\Rightarrow$ (3)]

 Suppose~$z$ is weakly 2-random and~$f$ is an a.e.\ differentiable computable function.  

Since $z$ cannot be contained in any $\arsigma03$ set of measure 0, it must belong to all $\arpi03$ sets of full measure. In particular, by Lemma \ref{pi03_lemma}, $z$ belongs to the set of differentiable points of $f$.

\end{proof}

\begin{proof}[Proof (3) $\Rightarrow$ (2)] \renewcommand{\qedsymbol}{} Trivial. \end{proof}

\begin{proof}[Proof (2) $\Rightarrow$ (1)] 

Suppose $z=(z_1,\dots,z_n)\in\zn$ is not weakly 2-random. 

We may assume that all coordinates of $z$ are weakly 2-random, otherwise the required conclusion follows  from the one dimensional case. For   suppose some $z_j$ is not weakly 2-random. Then there is a computable a.e.\ differentiable function $g:\zeroone\to\real$ such that $g'(z_j)$ does not exist. Define $\gamma:\zn\to \real$ as $\gamma(x_1,\dots,x_j,\dots x_n)=g(x_j)$. Then $\gamma$ is a computable a.e.\ differentiable function such that $\gamma'(z)$ doesn't exist.

In what follows, we ignore those elements of $\zn$ that have at least one of its coordinates rational. 

Let $\seq G$ be a sequence of uniformly $\arsigma01$ subsets of $\zeroone^n$ such that  $G_{i+1}\subseteq G_{i}$ for all $i$ and $G=\bigcap G_i$ is a null-set with $z\in G$. Since we ignore elements with dyadic coordinates, we may assume that every $G_i$ is an infinite union of basic dyadic $n-$cubes.

Let $\dseq D m l$ be an effective double sequence of (open) basic dyadic $n-$cubes such that $G_m=\bigcup_i D_{m,i}$ for each $m$, and  for all $n,k$ there is  an $l$ with \mbox{$D_{n+1,k}\subseteq D_{n,l}$}.

\vspace{5pt}
\n\emph{General idea of the proof.} We will construct a computable double sequence $\dseq C m i$ of basic dyadic $n-$cubes with certain well-behaved properties. For every $n-$cube $C_{m,i}$ in the sequence we will define a \emph{tent function} $f_{m,i}$ which is  0 outside $C_{m,i}$ and its graph forms a piecewise linear ``tent'' at $C_{m,i}$. See figure \ref{fig_tent} for an illustration of what a graph of a tent function on $\zeroone^2$ might look like. 

$E_{m,i}$ is the subarea of $C_{m,i}$ where $\left|D_1f_{m,i}\right|\neq \pm1$. The tent functions are defined in such a way that $z$ belongs to only finitely many of $E_{m,i}$.   This is where our assumption that all coordinates of $z$ are not weakly 2-random is used. See figure \ref{fig2}. Any point belonging to infinitely many $E_{m,i}$, by pigeonhole principle, must have at least one coordinate belonging to the (darker) corner areas (one-dimensional $E^k_{m,i}$ sets). Our tent functions are defined in such a way that those areas form $\arpi02$ null sets.

 Then $f:\zn\to\real$ will be defined as a sum of those $f_{m,i}$ for which we know the first partial derivative on $z$ is equal to $\pm1$. This is used to show that $D_1f(z)$ does not exist. The properties of $\dseq C m i$  ensure that $f$ is computable and a.e.\ differentiable.

\definecolor{ttttqq}{rgb}{0.62,0.62,0.62}
\definecolor{zzttqq}{rgb}{0.71,0.72,0.71}
\definecolor{xdxdff}{rgb}{0.049019607843137253,0.049019607843137253,0}
\definecolor{qqqqff}{rgb}{0.0,0.0,0.2}

\begin{figure}

\begin{tikzpicture}[line cap=round,line join=round,>=triangle 45,x=1.0cm,y=1.0cm]

\clip(0.0,1.0) rectangle (12.0,7.0);
\fill[color=zzttqq,fill=zzttqq,fill opacity=0.1] (1.84,1.44) -- (5.41,4.91) -- (11.41,4.91) -- (7.84,1.44) -- cycle;
\fill[color=zzttqq,fill=zzttqq,fill opacity=0.1] (2.97,6.62) -- (1.84,1.44) -- (5.41,4.91) -- cycle;
\fill[color=zzttqq,fill=zzttqq,fill opacity=0.1] (6.71,6.62) -- (7.84,1.44) -- (11.41,4.91) -- cycle;
\fill[color=zzttqq,fill=zzttqq,fill opacity=0.1] (2.97,6.62) -- (6.71,6.62) -- (7.84,1.44) -- (1.84,1.44) -- cycle;
\fill[line width=2.4000000000000004pt,color=ttttqq,fill=ttttqq,fill opacity=0.15] (2.97,1.44) -- (6.54,4.91) -- (10.28,4.91) -- (6.71,1.44) -- cycle;

\draw [color=black, dashed] (1.84,1.44)-- (5.41,4.91);
\draw [color=black, dashed] (5.41,4.91)-- (11.41,4.91);
\draw [color=black, dashed] (11.41,4.91)-- (7.84,1.44);
\draw [color=black, dashed] (7.84,1.44)-- (1.84,1.44);
\draw [color=black, dashed] (2.97,6.62)-- (1.84,1.44);
\draw [color=black, dashed] (1.84,1.44)-- (5.41,4.91);
\draw [color=black, dashed] (5.41,4.91)-- (2.97,6.62);
\draw [color=black] (6.71,6.62)-- (7.84,1.44);
\draw [color=black] (7.84,1.44)-- (11.41,4.91);
\draw [color=black] (11.41,4.91)-- (6.71,6.62);
\draw [color=black] (2.97,6.62)-- (6.71,6.62);
\draw [color=black] (6.71,6.62)-- (7.84,1.44);
\draw [color=black] (7.84,1.44)-- (1.84,1.44);
\draw [color=black] (1.84,1.44)-- (2.97,6.62);
\draw [line width=1.4000000000000004pt,color=black, dashed] (2.97,1.44)-- (6.54,4.91);
\draw [line width=1.4000000000000004pt,color=black, dashed] (6.54,4.91)-- (10.28,4.91);
\draw [line width=1.4000000000000004pt,color=black, dashed] (10.28,4.91)-- (6.71,1.44);
\draw [line width=1.4000000000000004pt,color=black] (6.71,1.44)-- (2.97,1.44);
\draw (2.97,6.62)-- (2.97,1.44);
\draw (6.71,6.62)-- (6.71,1.44);
\draw (1.84,1.44)-- (2.97,1.44);
\draw (6.71,1.44)-- (7.84,1.44);

\begin{scriptsize}
\begin{scope}[ opacity                 = 1.0
                 , execute at begin node = $\displaystyle 
                 , execute at end node   = $]

\draw [fill=qqqqff] (7.84,1.44) circle (1.5pt);
\draw[color=black] (7.951766200060514,1.672598071664566) node[right] {(a_1,b_2)};
\draw [fill=qqqqff] (1.84,1.44) circle (1.5pt);
\draw[color=black] (1.9592008926267068,1.672598071664566) node[left] {(a_1, a_2)};
\draw [fill=qqqqff] (5.41,4.91) circle (1.5pt);
\draw[color=black] (5.5313707416262465,5.144613625832338) node {(b_1,a_2)};
\draw [fill=qqqqff] (11.41,4.91) circle (1.5pt);
\draw[color=black] (11.523936049060055,5.144613625832338) node {(b_1,b_2)};
\draw [fill=xdxdff] (2.97,1.44) circle (1.5pt);
\draw[color=xdxdff] (3.0942829007200183,1.672598071664566) node {~};
\draw [fill=xdxdff] (6.54,4.91) circle (1.5pt);
\draw[color=xdxdff] (6.649760367247597,5.144613625832338) node {~};
\draw [fill=qqqqff] (6.71,1.44) circle (1.5pt);
\draw[color=qqqqff] (6.833376574439162,1.672598071664566) node {~};
\draw [fill=qqqqff] (10.28,4.91) circle (1.5pt);
\draw[color=qqqqff] (10.388854040966743,5.144613625832338) node {~};

\draw[text=black] (6.866761339383084,3.325143936388649) node {C_{m,i}\setminus E_{m,i}};

\draw[text=black] (2.44327998431356,1.322058039753397) node {\epsilon};
\draw[text=black] (7.317455666126016,1.322058039753397) node {\epsilon};
\end{scope}

\end{scriptsize}

\end{tikzpicture}

\caption{Two-dimensional graph of a tent function $f_{m,i}$.}\label{fig_tent}
\end{figure}
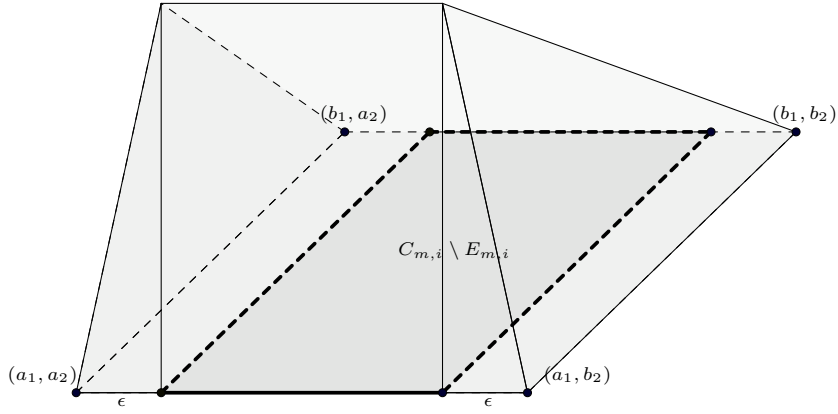

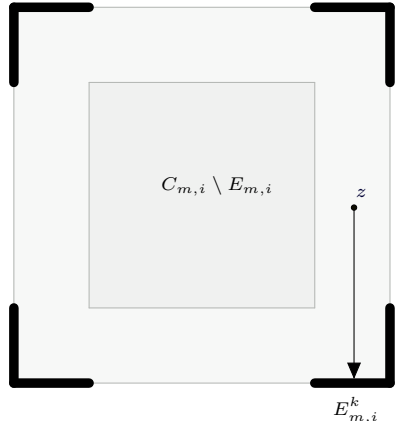
\begin{figure}[htb]
\begin{center}

\begin{tikzpicture}[line cap=round,line join=round,>=triangle 45,x=1.0cm,y=1.0cm]
\clip(0.5,-0.5) rectangle (6.5,5.5);
\fill[color=zzttqq,fill=zzttqq,fill opacity=0.1] (1.0,5.0) -- (1.0,0.0) -- (6.0,0.0) -- (6.0,5.0) -- cycle;
\fill[color=zzttqq,fill=zzttqq,fill opacity=0.1] (2.0,4.0) -- (2.0,1.0) -- (5.0,1) -- (5.0,4.0) -- cycle;
\draw [color=zzttqq] (1.0,5.0)-- (1.0,0.0);
\draw [color=zzttqq] (1.0,0.0)-- (6.0,0.0);
\draw [line width=0.4pt,color=zzttqq] (6.0,0.0)-- (6.0,5.0);
\draw [line width=0.4pt,color=zzttqq] (6.0,5.0)-- (1.0,5.0);
\draw [color=zzttqq] (2.0,4.0)-- (2.0,1.0);
\draw [color=zzttqq] (2.0,1.0)-- (5.0,1);
\draw [color=zzttqq] (5.0,1)-- (5.0,4.0);
\draw [color=zzttqq] (5.0,4.0)-- (2.0,4.0);
\draw [line width=3.6pt] (5.0,5.0)-- (6.0,5.0);
\draw [line width=3.6pt] (5.0,0.0)-- (6.0,0.0);
\draw [line width=3.6pt] (2.0,0.00)-- (1.0,0.0);
\draw [line width=3.6pt] (1.0,1.0)-- (1.0,0.0);
\draw [line width=3.6pt] (1.0,4.0)-- (1.0,5.0);
\draw [line width=3.6pt] (2.0,5.00)-- (1.0,5.0);
\draw [line width=3.6pt] (6.0,4)-- (6.0,5.0);
\draw [line width=3.6pt] (6.0,1)-- (6.0,0.0);
\draw [->] (5.522103681442526,2.3349538965917604) -- (5.522103681442526,0.05368895567242715);
\begin{scriptsize}
\begin{scope}[ opacity                 = 1.0
                 , execute at begin node = $\displaystyle 
                 , execute at end node   = $]

\draw[text=black] (3.7052878901714386,2.6218195478450896) node {C_{m,i}\setminus E_{m,i}};
\draw[color=black] (5.549424219657128,-0.08291373540058666) node[below] {E^k_{m,i}};
\draw [fill=xdxdff] (5.522103681442526,2.3349538965917604) circle (1pt);
\draw[color=qqqqff] (5.617725565193636,2.52619766409398) node {z};
\end{scope}
\end{scriptsize}
\end{tikzpicture}
\end{center}
\caption{Two-dimensional projection of $C_{m,i}$}
\label{fig2}
\end{figure}

\vspace{5pt}
\n\emph{Construction of the double sequence $\dseq C m i$.}

Suppose $m=0$, or $m>0$ and we have already defined $\left(C_{m-1,j}\right)_{j\in \NN}$. Define $\left(C_{m,j}\right)_{j\in \NN}$ as follows.

Let $N\in\NN$ be  the greatest number such that we have already defined $C_{m,i}$ for $i\le N$. When a new $n-$cube $D=D_{m,l}$ is enumerated into $G_m$, if $m>0,$ we wait until $D$ is contained in a union of $n-$cubes $\bigcup_{r\in F} C_{m-1,r}$, where $F$ is finite. This is possible since $D$ is contained in a single cell of the form $D_{m-1,\_}$, that was handled in a previous stage. If $m=0$, let $\delta=\lm{D}$, otherwise let $$\delta=\min\{\lm{D},\min\{\lm{C_{m-1,r}}:r\in F\}\}.$$  If $N=0$, let $\epsilon=8^{-m}\delta$, otherwise let $\epsilon=\min\{8^{-m}\delta,\lm{C_{m,N-1}}\}$. 

Finally, partition $D$ into disjoint basic dyadic $n-$cubes $C_{m,i}$, $i=N+1,\dots,N'<~\infty$, with nonincreasing volume $\lambda(C_{m,i})\le\epsilon^n$, so that when $m>0$, each of the cubes is contained in one of $C_{m-1,r}$ for some $r\in F$. 

The following claim summarizes all properties of $\dseq C m i$ relevant to our proof.
\begin{claim}
\n The double sequence $\dseq C m i$ is computable and it verifies the following properties:
\begin{enumerate}
\item[i)]  $G_m=\bigcup_{i\in \NN} C_{m,i},$ 
\item[ii)]$C_{m,i}~\cap C_{m,k}=\emptyset\text{ and } \lambda(C_{m,i})\ge\lambda(C_{m,k})\text{ for all }i<k,$ 
\item[iii)] if $B=C_{m,i}$ for $m>0$, then there is an $n-$cube $A=C_{m-1,k}$ such that 
\begin{align}\label{l17}
B\subseteq A\text{ and } \lambda(B)\le 8^{-m}\lambda(A), 
\end{align}
\item[iv)] for all $m,k\in\NN$ $$D_{m,k}= \text{some finite union of n-cubes of the form } C_{m,i}$$ with 
\begin{align}
C_{m,i}\subseteq D_{m,k}\implies d_{m,i}\le 8^{-m}\lm{D_{m,k}}
\end{align}
where $d_{m,i}$ denotes the length of a side of $C_{m,i}$.
\end{enumerate} 

\begin{proof} All of the listed properties are straightforward consequences of the construction of $\dseq C mi$.
\end{proof}
\end{claim}

\vspace{5pt}
\n \emph{Tent functions $f_{m,j}$.}

Let $m,j\in\NN$. For all $i\in\N$ with $1\le i\le n$, define $a^i_{m,j},b^i_{m,j}$ so that $\left(a^i_{m,j},b^i_{m,j}\right)=~\pi_i(C_{m,j})$, where $\pi_i:\R^n\to\R$ denotes the projection onto the $i-$th coordinate.

Let $\epsilon_{m,j}=\epsilon=2^{-m-j-1}\cdot d_{m,j}$ and define $b^i_{m,j}:\zeroone\to\real$ as  
\begin{align*}
b^i_{m,j}(x) = \begin{cases} \frac{x-a^i_{m,j}}{\epsilon} &\mbox{if } x\in [a^i_{m,j},a^i_{m,j}+\epsilon], \\ 
1 & \mbox{if } x\in (a^i_{m,j}+\epsilon,b^i_{m,j}-\epsilon),\\ 
\frac{b^i_{m,j}-x}{\epsilon} &\mbox{if } x\in [b^i_{m,j}-\epsilon,b^i_{m,j}],\\
0&\mbox{otherwise}.
\end{cases} 
\end{align*}

Define $f_{m,j}:\zn\to\real$ as 
\begin{align*}
f_{m,j}(x_1,x_2,\dots,x_n) = d\left(\zeroone\setminus \left(a^1_{m,j},b^1_{m,j}\right),x_1\right) \cdot  \prod_{n\ge i\ge 2} b^i_{m,j}(x_i),
\end{align*}

where $ d\left(\zeroone\setminus \left(a^1_{m,j},b^1_{m,j}\right),x_1\right)$ denotes the distance from $x_1$ to $\zeroone\setminus \left(a^1_{m,j},b^1_{m,j}\right)$.

Note that $f_{m,j}$ is a computable (uniformly in $m,j$) a.e.\ differentiable function. 

Lastly, define $E_{m,j}$ to be the subset of $C_{m,j}$ where $\left|D_1f_{m,j}\right|\neq 1$ whenever $D_1f_{m,j}$ exists, that is $$E_{m,j}=C_{m,j} \setminus \left[\left(a^1_{m,j},b^1_{m,j}\right)\times \prod_{n\ge i\ge 2} \left(a^i_{m,j}+\epsilon, b^i_{m,j}-\epsilon\right) \right].$$

 The idea behind such definition of $f_{m,j}$ functions is that $\epsilon_{m,j}$ goes to $0$ so quickly, that $\left|D_1f_{m,j}(z)\right|\neq 1$ holds only for finitely many $m,j\in \NN$.

\begin{claim}
There exists $N\in\NN$ such that for all $i\in\NN$ and $m>N,$ if $z\in C_{m,i}$ then $\left|D_1f_{m,i}(z)\right|= 1$.
\begin{proof} 

To prove this claim, we will use our assumption that all coordinates of $z$ are weakly $2-$random. Specifically, we will show that if a point belongs to an infinitely  many $E_{m,i}$, then one of its coordinates belongs to a $\arpi02$ null set.

\n For every $m,i,k\in\NN$ with $2\le k\le n$, let $$E^k_{m,i}=\left(a^k_{m,i},a^k_{m,i}+\epsilon_{m,i}\right)\cup \left(b^k_{m,i}-\epsilon_{m,i}, b^k_{m,i}\right).$$ Note the following property of those sets: if $z\in E_{m,i}$ then for some $k$, $z_k\in E^k_{m,i}$.

For every $m,k\in\NN$ with $n\ge k\ge 2$, let $B_m^k=\bigcup_{i>m} \bigcup_j E^k_{i,j}$. Let's verify that every $B^k=\bigcap_i B_i^k$ is a $\arpi02$ null-set. Indeed, $\seq {B^k}$ is a uniformly computable sequence of $\arsigma 01$ sets with $\lm{B_m^k}\le \sum_{i>m}\sum_j 2^{-i-j}\cdot d_{i,j}\le 8^{-m}$ for all $m,k$.

By the pigeonhole principe, if $z$ belongs to infinitely many $E_{m,j}$ (for infinitely many $m$), then for some $k$, $z_k$ belongs to infinitely many $E^k_{m,j}$. In that case $z_k\in B^k$ and we get a contradiction.

 Let $N$ be such that $z_k\notin B_{N}^k$ for all $k$ and the required result follows.
\end{proof}
\end{claim}

\vspace{5pt}
\n \emph{Definition of the function $f$.} 

\n  Let $$f_m=\sum_{i=0}^\infty 4^{m}f_{m,i}$$ and $$f=\sum_{i> N} f_m.$$

\vspace{5pt}
\n

\begin{claim} $f$ is computable.
\begin{proof}  
 
Fix $m>0$. Note that every $f_{m,i}$ is bounded from above by $d_{m,i}/2$ and since all $C_{m,i}$ are disjoint, $f_m$ is bounded from above by $4^m8^{-m}/2=2^{-m-1}$ and it follows that $f$ is well defined everywhere.

Firstly, let's show that $f(q)$ is computable uniformly in rational $q$. Given $m>0$, since $\lim_{i\to\infty}\lm{C_{m,i}}=0,$ we can find $i^*$ such that $${d_{k,i^*}}\le 8^{-m}/(m+1)\text{ for each $k\le m$}.$$
Since  the $d_{k,i}$ is non-increasing in $i$ and $ {f_{k,i}}\le d_{k,i}/2$, we have
$$4^k f_{k,i} (q)\le 2^{-m-1}/(m+1)\text{ for all $k\le m$ and $i\ge i^*$}.$$ Hence  $$\sum_{k\le m}\sum_{i\ge i^*}4^k {f_{k,i}}(q)\le 2^{-m-1}.$$    Furthermore,  $$\sum_{k>m}f_k(q)\le \sum_{k>m}2^{-k-1}=2^{-m-1}.$$ Therefore the approximation of $f(q)$ at stage $i^*$ based only on the $n-$cubes of the form $C_{k,i}$ for $k\le m$ and $i< i^*$ is within $2^{-m}$ of $f(q)$.

\

\n Secondly, we need to verify that $f$ is effectively uniformly continuous. Suppose $\|x-y\|\le d_{m,1}$ for some $m$. Then for $k<m$, we have $|f_k(x)-f_k(y)|\le 4^{k}d_{m,1}/2.$ For $k\ge m$, we have  $f_k(x),f_k(y)\le 2^{-k-1}$.  Thus $$|f(x)-f(y)|\le d_{m,1}\sum_{k<m}4^k +\sum_{k\ge m}2^{-k}<2^{-m+2}.$$

\n Define $h(m)=\lfloor-\log_2d_{m,1}\rfloor+1$ so that $2^{-h(m)}\le d_{m,1}$. Note that $h$ is a computable order function.  Then we get that $\|x-y\|\le 2^{-h(m)}$ implies $|f(x)-f(y)|\le 2^{-m+2}.$
 
\end{proof}
\end{claim}

\begin{claim}\label{claim133} $D_1f(z)$ does not exist.
\begin{proof}

    For all $m>N,$ let $d_m=d_{m,i_m}$ where $i_m$ is such that $z\in C_{m,i_m}$. Note that for all $m>N$ we have either  $\delta_{f_m}^1\left(z, \frac{{d_m}}{4}\right) = \pm 4^m$  or $\delta_{f_m}^1\left(z, -\frac{{d_m}}{4}\right)= \pm 4^m$. Without loss of generality we may assume that for infinitely many $m$, we have $\left| \delta_{f_m}^1\left(z, \frac{{d_m}}{4}\right)\right|=  4^m.$ Fix one such $m>N$.
Note that for for every $k\in\NN$ with $N<k<m$ we have $\left|\delta_{f_k}^1\left(z, \frac{{d_m}}{4}\right)\right|= 4^k$. Suppose $k>m$. Then we have $$f_k(x)\le 4^k8^{-k}\frac{d_m}{2}=2^{-k-1}d_m$$ for all $x\in C_{m,l}\setminus E_{m,l} $ and thus we get

$$\left|\delta_{f_k}^1\left(z, \frac{{d_m}}{4}\right)\right|\le \frac{2\cdot2^{-k-1} d_m}{\left\|\frac{{d_m}}{4} e_1\right\|}=2^{-k+2}.$$ 

\n Hence, for $m>N$ we have $$\left|\delta_{f}^1\left(z, \frac{{d_m}}{4}\right)\right|\ge \left(4^m-\sum_{N<k<m}4^k-\sum_{k>m}2^{-k+2}\right)\ge 4^{m-1}-4.$$
Therefore $D_1f(z)$ does not exist.

\end{proof}
\end{claim}

\begin{claim}
$f$ is differentiable almost everywhere.
\begin{proof} 
Let $x\in\zn$. There are three possible cases:
\begin{enumerate}
\item $f_{m,j}'(x)$ does not exist for some $m,j$,
\item $x$ belongs to the support of $f_{m,j}$ for infinitely many $m,j$, or
\item  $x$ belongs to the support of $f_{m,j}$ for only finitely many $m,j$ and all $f_{m,j}'(x)$ exist. Note that this implies differentiability of $f$ at $x$.
\end{enumerate}

The first case corresponds to a null-set, since every $f_{m,j}$ is a.e.\ differentiable. The second case corresponds to a null-set too, since it implies $x\in \bigcap_i G_i$. The last case implies differentiability of $f$ at $x$ and it must correspond to a set   of full measure since the cases (1) and (2) are captured by null-sets. Thus $f$ is a.e.\ differentiable. 
\end{proof}
\end{claim}

\end{proof}

\end{theorem}

\section{Conclusion and future directions}

Despite the obstacle described in the Subsection \ref{dore_maleva}, we conjecture that computable randomness, just like on the unit interval,  characterises differentiability points of all computable real-valued Lipschitz functions.
Proving the converse to the effective version of Rademacher's Theorem (that is, showing that differentiability of computable Lipschitz functions implies computable randomness) remains an open question of great interest.

There are quite a few  results in classical analysis about differentiability of functions of several variables that exhibit Lipschitz-like behaviour. Naturally, those results are related to Rademacher's Theorem. Studying effective versions of those will improve our understanding of interplay between computable analysis and algorithmic randomness. Here we mention two such theorems that we feel are of particular importance: 

\n(1) Alexandrov's theorem (see 6.4 in \cite{Evans.Gariepy:92}) states that convex functions are twice differentiable almost everywhere. Convex functions and monotone functions are closely related: on the real line, a function is monotone if and only if it is a derivative of a convex function. For functions of several variables, the relation is a bit less straightforward    (see \cite{Rockafellar:70}). Recently it has been shown that both twice-differentiability of computable convex real valued functions on $\R^n$ and differentiability computable monotone functions on $\R^n$ correspond to computable randomness (see \cite{Galicki:15,Galicki:15:2}).

\n(2) It is known that $K_n$-monotone functions of several variables are a.e.\ differentiable (see \cite{Chabrillac.Crouzeix:87}). Two of the three steps of our proof in Section \ref{lip_section} work for $K_n-$monotone functions. The one that doesn't work is the one in Subsection \ref{existence_directional_derivatives}. It is not known whether computable randomness implies differentiability of $K_n-$monotone computable functions, and what randomness notion is induced by a.e. differentiability of computable $K_n-$monotone functions.

On the other hand, our result concerning weak 2-randomness is sharp: weak 2-randomness does characterise differentiability sets of computable a.e.\ differentiable functions of several variables. There are many other similar results in one dimension that characterise differentiability of effective functions in terns of algorithmic randomness. Generalizing those results to higher dimensions (and, perhaps, to more general spaces) will provide more insight into interactions between computable analysis and algorithmic randomness. 

\section*{Acknowledgements}

The first author  would like to thank Andr\'e Nies for discussing the paper and for many helpful suggestions and corrections, and Kenshi Miyabe and Jason Rute for explaining their  results on computable randomness and for valuable discussions.

Part of this work was done while both co-authors were invited to the Institute for Mathematical Sciences of the National University of Singapore in June 2014.

\bibliographystyle{plain}

\end{document}